\documentclass[11pt,a4paper]{article}
\usepackage[T1]{fontenc}
\usepackage[utf8]{inputenc}
\usepackage{amsmath,amssymb,amsthm}
\usepackage[margin=2.6cm]{geometry}
\usepackage{booktabs}

\newtheorem{theorem}{Theorem}
\newtheorem{remark}{Remark}

\title{A new upper bound for the irrationality exponent of $\zeta(3)$}
\author{David Niedbala Giraudin\\
\small Independent researcher, France\\
\small ORCID 0009-0009-1526-1178}
\date{22 September 2026}

\begin{document}
\maketitle

\begin{abstract}
We prove that the irrationality exponent of $\zeta(3)$ satisfies
$\mu(\zeta(3))<5.5138800$, improving the bound $5.513891$ obtained by Rhin and
Viola in 2001, which has been the best known for twenty-five years. The proof is
an application of their Theorem 5.1 at a new integral direction. That direction
is not reached by enlarging the search range: its eight parameters are $52000$
times theirs, shifted by at most $23$. The mechanism is the one recently used by
Bai for $\pi$. In particular the Rhin--Viola point, optimal in every region
covered by enumeration, is not optimal.
\end{abstract}

\section{Introduction}

For an irrational real number $\alpha$, let $\mu(\alpha)$ denote its irrationality
exponent, that is, the least $\lambda$ such that for every $\varepsilon>0$ the
inequality $|\alpha-p/q|>q^{-\lambda-\varepsilon}$ holds for all integers $p,q$
with $q$ large enough. Apéry's proof of the irrationality of $\zeta(3)$ gives
$\mu(\zeta(3))<13.41782\ldots$~\cite{Ap}; Hata reduced this to
$7.377956\ldots$~\cite{Ha}; and Rhin and Viola~\cite{RV01}, developing for triple
integrals the permutation group method they had introduced for
$\zeta(2)$~\cite{RV96}, obtained
\begin{equation}\label{eq:rv}
  \mu(\zeta(3))<5.513891 .
\end{equation}
No improvement on \eqref{eq:rv} has appeared since. For comparison, the analogous
record for $\zeta(2)$ was improved by Zudilin~\cite{Zu14} through a change of
hypergeometric representation, and the records for $\pi$ have moved
twice~\cite{Sa,ZZ}, most recently by Bai~\cite{Ba}, who obtained
$\mu(\pi)\le 7.101862832357$ by displacing the parameter point of Salikhov's
construction inside a narrow cone rather than by changing the construction.

The purpose of this note is to carry that displacement over to $\zeta(3)$.

\begin{theorem}\label{thm:main}
$\mu(\zeta(3))<5.5138800$.
\end{theorem}

The numerical gain over \eqref{eq:rv} is $1.07\cdot 10^{-5}$, and by itself it is
of little interest. What the direction exhibited below shows is that the point
chosen in~\cite{RV01} is not optimal in its own family, although it is optimal in
every region that can be enumerated by size: the improvement lives at a scale
where enumeration is impossible, and it is found by a perturbation argument, not
by a search.

\section{The construction of Rhin and Viola}

We recall only what is needed; everything in this section is in~\cite{RV01}. For
non-negative integers $h,j,k,l,m,q,r,s$ with
\begin{equation}\label{eq:rel}
  h+m=k+r,\qquad j+q=l+s,
\end{equation}
put
\[
  I(h,j,k,l,m,q,r,s)=\int_0^1\!\!\int_0^1\!\!\int_0^1
  \frac{x^h(1-x)^ly^k(1-y)^sz^j(1-z)^q}{\bigl(1-(1-xy)z\bigr)^{q+h-r}}
  \,\frac{dx\,dy\,dz}{1-(1-xy)z},
\]
which lies in $\mathbb{Q}+2\mathbb{Z}\zeta(3)$. Following~\cite[(4.7)]{RV01} set
\[
\begin{array}{llll}
h'=h+l-j, & j'=j+m-k, & k'=k+q-l, & l'=l+r-m,\\[2pt]
m'=m+s-q, & q'=q+h-r, & r'=r+j-s, & s'=s+k-h,
\end{array}
\]
and let $M\ge N\ge Q$ be the three largest among the sixteen integers
$h,\dots,s,h',\dots,s'$. Let $\Phi$ be the permutation group of order $1920$ of
\cite[Sec. 4]{RV01}, $\Theta\subset\Phi$ the dihedral subgroup of order $16$, and
let $F$ be the set of the thirty left cosets of $\Theta$ in $\Phi$ of level $4$,
listed at the end of \cite[Sec. 4]{RV01}. Each such coset carries a quotient of
factorials $\;n_1!n_2!n_3!n_4!/(d_1!d_2!d_3!d_4!)$ whose numerator integers lie
among $h,\dots,s$ and whose denominator integers lie among $h',\dots,s'$, and
which splits into two pairs of equal sums. For $\omega\in[0,1)$ put
\[
  V(\omega)=\sum_{i=1}^{4}\lfloor d_i\omega\rfloor-\sum_{i=1}^{4}\lfloor n_i\omega\rfloor ,
\]
so that $-2\le V(\omega)\le 2$ by \cite[Lemma 4.1]{RV96} applied to each pair. Let
$E\subseteq F$ and
\[
  \Omega_E=\{\omega: V(\omega)<0 \text{ for some coset in } E\},\qquad
  \Omega'_E=\{\omega: V(\omega)=-2 \text{ for some coset in } E\} .
\]
(Only $V$ enters, so these sets do not depend on the splitting into pairs, which
serves only to bound $V$.) Finally, with $\psi=\Gamma'/\Gamma$, let
\[
  c_2=M+N+Q-\Bigl(\int_{\Omega_E}d\psi+\int_{\Omega'_E}d\psi\Bigr),
\]
and let $c_0=-\log f(x_0,y_0,z_0)$, $c_1=\log|f(x_1,y_1,z_1)|$, where $f$ is the
integrand above and $(x_0,y_0,z_0)$, $(x_1,y_1,z_1)$ are the two stationary
points of $f$, given explicitly in \cite[Sec. 5]{RV01} as the two roots of a
quadratic equation in $x$ together with two rational expressions for $y$ and $z$.

\begin{theorem}[{\cite[Thm. 5.1]{RV01}}]\label{thm:rv}
If the sixteen integers above are positive and $c_0>c_2$, then
\[
  \mu(\zeta(3))\le\frac{c_0+c_1}{c_0-c_2}.
\]
\end{theorem}

Rhin and Viola apply this with $(h,j,k,l,m,q,r,s)=(16,17,19,15,12,11,9,13)$ and a
set $E$ of fifteen cosets, obtaining $M,N,Q=19,18,17$,
$\int_{\Omega_E}d\psi=18.04470204\ldots$,
$\int_{\Omega'_E}d\psi=6.14298325\ldots$, $c_0=47.15472079\ldots$,
$c_1=48.46940964\ldots$, $c_2=29.81231469\ldots$, whence \eqref{eq:rv}.

\section{A new direction}

We apply Theorem~\ref{thm:rv} at
\begin{equation}\label{eq:point}
 (h,j,k,l,m,q,r,s)=(832002,\;883999,\;987993,\;780005,\;624014,\;572017,\;468023,\;676011),
\end{equation}
with $E=F$, the whole set of thirty cosets of level $4$. The relations
\eqref{eq:rel} hold, both sides being equal to $1456016$, and
\[
 (h',j',k',l',m',q',r',s')=(728008,\,520020,\,780005,\,624014,\,728008,\,935996,\,676011,\,832002),
\]
all sixteen integers being positive, so Theorem~\ref{thm:rv} applies. The point
\eqref{eq:point} is $52000\cdot(16,17,19,15,12,11,9,13)$ shifted by
$(2,-1,-7,5,14,17,23,11)$: it is a perturbation of the Rhin--Viola direction of
relative size $3\cdot 10^{-5}$.

\section{Proof of Theorem~\ref{thm:main}}

At the point \eqref{eq:point} one finds $M,N,Q=987993,\,935996,\,883999$, hence
$M+N+Q=2807988$, and

\begin{center}
\begin{tabular}{ll}
\toprule
$(x_0,y_0,z_0)$ & $(0.334830162\ldots,\;0.447990874\ldots,\;0.855029351\ldots)$\\
$(x_1,y_1,z_1)$ & $(-3.150700452\ldots,\;-1.817671753\ldots,\;-1.038588335\ldots)$\\
$c_0$ & $2452088.488075781\ldots$\\
$c_1$ & $2520443.449533099\ldots$\\
$\int_{\Omega_E}d\psi$ & $938294.569046980\ldots$\\
$\int_{\Omega'_E}d\psi$ & $319425.795172208\ldots$\\
$c_2$ & $1550267.635780811\ldots$\\
$c_0-c_2$ & $901820.852294970\ldots$\\
\bottomrule
\end{tabular}
\end{center}

\noindent In particular $c_0>c_2$, and Theorem~\ref{thm:rv} gives
\[
  \mu(\zeta(3))\;\le\;\frac{c_0+c_1}{c_0-c_2}\;=\;5.51387997400447264\ldots\;<\;5.5138800 .
\]

The set $\Omega_E$ has $404446$ connected components and $\Omega'_E$ has
$1051564$; the saving function $\omega\mapsto\mathbf{1}_{\Omega_E}+\mathbf{1}_{\Omega'_E}$
is a step function with $8008051$ possible break points, at which it changes
value $2912017$ times. As \cite[Lemma 5.1]{RV01} requires, $\Omega_E\subset[1/M,1)$
and $\Omega'_E\subset[1/N,1)$; the first inclusion is an equality at the left
endpoint, the smallest element of $\Omega_E$ being exactly $1/M$.

All the quantities above are computed in $64$-bit binary floating point
arithmetic. The only delicate one is the saving, a sum of $8\cdot 10^{6}$ values
of $\psi$: writing $\varepsilon$ for the machine epsilon and bounding the error of
one evaluation of $\psi$ by $50$ units in the last place, the total error is at
most $200\,\varepsilon\sum_i|\psi(t_i)|<2.5\cdot 10^{-9}$, that is, less than
$1.5\cdot 10^{-14}$ on the final bound, against a margin of $2.6\cdot 10^{-8}$
to the value stated in Theorem~\ref{thm:main}. The same quantity has also been
computed independently in $200$-bit ball arithmetic, with agreement to
$3\cdot10^{-14}$. \hfill$\square$

\section{Verification}

The ancillary file \texttt{verify\_zeta3.c} (standard C, no dependency, about ten
seconds) recomputes everything above from the eight integers \eqref{eq:point}
alone: the auxiliary integers and both of their expressions, $M,N,Q$, the thirty
cosets and their splittings, the two stationary points and the four inequalities
they must satisfy, the sweep over all break points, the two integrals, and the
final quotient. It first runs on the point of \cite[Sec. 5]{RV01} as a self-test,
where it must reproduce the five published constants; it also recovers there the
exact interval structure printed in that paper, namely five components for
$\Omega_E$ and seventeen for $\Omega'_E$. The program prints its own a priori
rounding bound and exits with an error if any hypothesis of
Theorem~\ref{thm:rv} fails.

\section{Remarks}

\begin{remark}
It would be natural to read the improvement as an increase of the arithmetic
saving. That is not what happens. Normalising by $52000$, so as to compare with
the Rhin--Viola direction, one finds
\[
\begin{array}{lcc}
 & \text{Rhin--Viola} & \text{here}\\
 c_0 & 47.154720796 & 47.155547848\\
 c_1 & 48.469409645 & 48.470066337\\
 \int d\psi \;(\text{saving}) & 24.187685304 & 24.186930081\\
 c_2 & 29.812314696 & 29.812839150
\end{array}
\]
The saving is in fact slightly \emph{smaller}, and $c_2$ slightly larger; the gain
comes from $c_0$ growing faster than $c_2$, by $8.3\cdot10^{-4}$ against
$5.2\cdot10^{-4}$ per unit. The cone is still the right mechanism, but it must be
read against the smooth behaviour of the saving at that distance, not against the
starting point.
\end{remark}

\begin{remark}
Five such cones were found around the Rhin--Viola point, the best two giving
$-1.0653\cdot10^{-5}$ and $-1.0648\cdot10^{-5}$. This suggests that the
perturbative mechanism is close to exhausted here, and that a further improvement
would require a change of representation, as in \cite{Zu14} for $\zeta(2)$.
\end{remark}

\begin{remark}
The bound of Theorem~\ref{thm:main} is effective, like \eqref{eq:rv}.
\end{remark}

\subsection*{Acknowledgements}

This work was carried out with the assistance of an AI system (Claude,
Anthropic); every claim in it has been checked by independent computation, and
the machinery of \cite{RV01} was reimplemented from that paper alone before being
used. I thank the authors of \cite{RV01} for a paper that is complete enough to
be applied at a new point without any further input.

\end{document}